\documentclass[12pt]{amsart}
\usepackage{amssymb,latexsym,comment,url}

\newcommand{\Q}{{\mathbb Q}}

\newcommand{\rank}{\operatorname{rank}}

\newcommand{\Magma}{{\sf MAGMA }}

\newenvironment{Proof}{\par\noindent{\sc Proof:}}%
                      {\hspace*{\fill}\nobreak$\Box$\par\medskip}
                       {\hspace*{\fill}\nobreak$\Box$\par\medskip}

\newtheorem{Proposition}{Proposition}[section]
\newtheorem{Theorem}[Proposition]{Theorem}

\theoremstyle{definition}

\renewcommand{\baselinestretch}{1.1}

\begin{document}

\title[Sequences of consecutive squares on quartic elliptic curves]%
{Sequences of consecutive squares on quartic elliptic curves}

\author[M. Kamel]%
{Mohamed~Kamel}
\address{Department of Mathematics, Faculty of Science, Cairo University, Giza, Egypt}
\email{mohgamal@sci.cu.edu.eg}

\author[M. Sadek]%
{Mohammad~Sadek}
\address{American University in Cairo, Mathematics and Actuarial Science Department, AUC Avenue, New Cairo, Egypt}
\email{mmsadek@aucegypt.edu}

\begin{abstract}
Let $C: y^2=ax^4+bx^2+c$, be an elliptic curve defined over $\Q$. A set of rational points $(x_i,y_i) \in C(\Q)$, $i=1,2,\cdots,$ is said to be a sequence of consecutive squares if $x_i= (u + i)^2$, $i=1,2,\cdots$, for some $u\in \Q$. Using ideas of Mestre, we construct infinitely many elliptic curves $C$ with sequences of consecutive squares of length at least $6$. It turns out that these 6 rational points are independent. We then strengthen this result by proving that for a fixed $6$-term sequence of consecutive squares, there are infinitely many elliptic curves $C$ with the latter sequence forming the $x$-coordinates of six rational points in $C(\Q)$.
\end{abstract}

\maketitle

\section{Introduction}
In \cite{Bremner}, Bremner started discussing the existence of long sequences on elliptic curves. He produced an infinite family of elliptic curves with arithmetic progression sequences of length 8. Several authors displayed infinite families of elliptic curves with long arithmetic progression sequences, see \cite{Campbell,Macleod,Ulas1}.

A geometric progression sequence is another type of sequence that has been studied on elliptic and hyperelliptic curves. Infinitely many (hyper)elliptic curves with $5$-term and $8$-term geometric progression sequences have been introduced in \cite{BremnerUlas} and \cite{AlaaSadek} respectively.

In \cite{SadekKamel}, sequences of consecutive squares on elliptic curves were studied. Infinitely many elliptic curves defined by equations of the form $E: y^2=ax^3+bx+c$, $a,b,c\in \Q$, with 5-term sequences of consecutive squares were presented. This was achieved by identifying these curves as rational points on an elliptic surface whose rank is positive.

In this note, we discuss sequences of consecutive squares on elliptic curves defined by the equation $y^2=ax^4+bx^2+c$, $a,b,c\in\Q$. We construct infinitely many such curves with $6$-term sequences of consecutive squares. More precisely, given a $6$-term sequence of consecutive squares, we prove the existence of an elliptic curve on which the latter sequence forms the $x$-coordinates of six rational points. For the construction, we use an idea due to Mestre. This sequence corresponds to six linearly independent rational points on the elliptic curve. In particular, we give an infinite family of elliptic curves with $2$-torsion points and rank $\ge 6$.

Finally, given a fixed $6$-term sequence of consecutive squares $(t+i)^2$, $i=0,\pm1,\pm2,3$, we find infinitely many elliptic curves of the form $y^2=ax^4+bx^2+c$ for which $(t+i)^2$ is an $x$-coordinate of a rational point. This is performed by realizing these elliptic curves as rational points on an elliptic surface of positive Mordell-Weil rank.

\section{First construction}
Let $C$ be an elliptic curve defined over a number field $K$ by $y^2= P(x)$ where $P\in K[x]$ is a polynomial of degree either $3$ or $4$. The sequence $(x_i,y_i)\in C(K)$ is said to be a {\em sequence of consecutive squares} on $C$ if there is a $u\in K$ such that $x_i=(u+i)^2$, $i=1,2,\ldots$. The authors proved in \cite{SadekKamel} that this sequence must be finite.

  Mestre, \cite{Mestre}, constructed elliptic curves with Mordell-Weil rank $\ge 11$, using the following idea: For any monic polynomial $P\in \mathbb{Q}(x)$ of degree $2n$ there exists a monic polynomial $Q \in \mathbb{Q}(x)$ of degree $n$ and $R \in \mathbb{Q}(x)$ of degree at most $n-1$ such that $P=Q^2-R $. If $x \in \mathbb{Q}$ is a root of $P$, then there is a rational point $(x,Q(x))$ on the algebraic curve $y^2=R(x)$.

 \begin{Theorem}
 \label{thm1}
 For any nontrivial sequence of consecutive squares $\left(t-\frac{5}{2}\right)^2,\left(t-\frac{3}{2}\right)^2,\\ \left(t-\frac{1}{2}\right)^2, \left(t+\frac{1}{2}\right)^2, \left(t+\frac{3}{2}\right)^2,\left(t+\frac{5}{2}\right)^2$, there is an elliptic curve described by $E_t:y^2=a(t) x^4+b(t)x^2+c(t)$, $a,b,c\in\Q(t)$, such that $(t+i)^2$, $i=\pm \frac{1}{2},\pm \frac{3}{2}, \pm \frac{5}{2}$, is the $x$-coordinate of a rational point in $E_t(\Q(t))$. In particular, there are infinitely many elliptic curves described by $y^2=a x^4+bx^2+c$ with $6$-term sequences of consecutive squares.
\end{Theorem}
\begin{Proof}
Consider the degree $12$ polynomial
{\footnotesize\[P(x)=\Big(x^2-\Big(t-\frac{5}{2}\Big)^4\Big)\Big(x^2-\Big(t-\frac{3}{2}\Big)^4\Big)\Big(x^2-\Big(t-\frac{1}{2}\Big)^4\Big)\Big(x^2-\Big(t+\frac{1}{2}\Big)^4\Big)
\Big(x^2-\Big(t+\frac{3}{2}\Big)^4\Big)\Big(x^2-\Big(t+\frac{5}{2}\Big)^4\Big).\]}
One may write
$P(x)=Q(x)^2-R(x)$, where
{\footnotesize\begin{eqnarray*}
Q(x)&=&x^6+\frac{1}{16} \Big(-48 t^4-840 t^2-707\Big) x^4+\frac{1}{256} \Big(768 t^8+8960 t^6-9184 t^4-322000 t^2+51331\Big) x^2\\&+&\frac{-4096 t^{12}+71680 t^{10}-1994496 t^8-50973440 t^6-251212528 t^4-260162280 t^2-50625}{4096} ,
  \end{eqnarray*}}
  {\footnotesize\begin{eqnarray*}
   R(x)&=&9t^2 \Big\{\frac{1}{64} \Big(5376 t^{10}+779520 t^8+11657184 t^6+57509200 t^4+95561365 t^2+36613360\Big)x^4\\&-&\frac{1}{512} \Big(86016 t^{14}+6113280 t^{12}+71158528 t^{10}+145053440 t^8-1767894864 t^6-8757574840 t^4\\&-&7679989163 t^2+1441328880\Big) x^2+\frac{1}{16384}\Big(336 t^6+11320 t^4+54229 t^2+56560\Big) \Big(4096 t^{12}\\&-&71680 t^{10}+1220352 t^8+24892160 t^6+126268912 t^4+129848040 t^2+50625\Big) \Big\}.
\end{eqnarray*}}
We consider the elliptic curve $ E_t: y^2=R(x):=a(t)x^4+b(t)x^2+c(t)$. By definition of $P(x)$, the latter curve possesses the $6$ rational points $\left((t+i)^2,Q((t+i)^2)\right)$, $i=\pm \frac{1}{2},\pm \frac{3}{2},\pm \frac{5}{2}$.
\end{Proof}

Using the rational transformation $x\mapsto x-(t-\frac{1}{2})^2$, the elliptic curve $E_t:y^2=R(x)=a(t)x^4+b(t)x^2+c(t)$ may be described by an equation of the form $y^2=A(t)\,x^4+B(t)\,x^3+C(t)\,x^2+D(t)\,x+E(t)^2.$
By virtue of \cite[Proposition 1.2.1]{Connell}, the curve $E_t$ is birationally equivalent, hence isomorphic, over $\mathbb{Q}(t)$ to the curve $E_t^*:T^2=S(S^2+\alpha(t)S+\beta(t))$ which has a nontrivial 2-torsion point $(0,0)$, where
{\footnotesize\begin{eqnarray*}
  \alpha(t) &=& \frac{9}{256} t^2 \Big(86016 t^{14}+6113280 t^{12}+71158528 t^{10}+145053440 t^8-1767894864 t^6-8757574840 t^4\\&-&7679989163 t^2+1441328880\Big) ,\\
  \beta(t) &=& -\frac{243 t^4}{1024}\left(4 t^2+17\right) \left(4 t^2+33\right) \left(4 t^2+97\right) \left(28 t^2+151\right) \left(4 t^3-48 t^2+t-68\right) \\&\times&\left(4 t^3-24 t^2+9 t-26\right)
   \left(4 t^3+24 t^2+9 t+26\right) \left(4 t^3+48 t^2+t+68\right) \left(20 t^3-24 t^2+125 t-10\right)\\&\times& \left(20 t^3+24 t^2+125 t+10\right).
   \end{eqnarray*}}
The rational points $(x_i,y_i)=((t+i-\frac{5}{2})^2,Q((t+i-\frac{5}{2})^2))$, $i=0,1,2,3,4,5$, on the curve $E_t$ correspond to $q_i=(S_i,T_i)$ on the curve $E_t^*$, where
{\footnotesize\begin{eqnarray*}
  S_i &=& -\frac{1}{64 \left(4 t^2-4 t+4 x_i+1\right)^2}\Big\{-1024 D(t) t^2+1024 D(t) t-1024 D(t) x-256 D(t)-2048 E^2(t)\\&-&2048 E(t) y_i+3096576 t^{19}-122978304 t^{18}+6193152 t^{17}
   x_i+687992832 t^{17}-239763456 t^{16} x_i\\&-&3489546240 t^{16}+3096576 t^{15} x_i^2+1134673920 t^{15} x_i+12386571264 t^{15}-116785152 t^{14} x_i^2\\&-&5784477696
   t^{14} x_i-36507285504 t^{14}+449777664 t^{13} x_i^2+18704996352 t^{13} x_i+84282105600 t^{13}\\&-&2413264896 t^{12} x_i^2-52863455232 t^{12} x_i-186072109056
   t^{12}+6826788864 t^{11} x_i^2+111024506880\\&\times& t^{11} x_i+294253225920 t^{11}-19001622528 t^{10} x_i^2-247903847424 t^{10} x-482333956992 t^{10}\\&+&34803933696 t^9
   x_i^2+312846477696 t^9 x_i+595208516688 t^9-84397584384 t^8 x_i^2-589845474432 t^8 x_i\\&-&578162334528 t^8+63324671040 t^7 x_i^2+522359939520 t^7 x_i+478308410748
   t^7-210498670080 t^6 x_i^2\\&-&486503360928 t^6 x_i-249506322696 t^6+34850131920 t^5 x_i^2+339523475688 t^5 x_i+37405300797 t^5\\&-&155776881024 t^4 x_i^2-37863329472
   t^4 x_i+25555102056 t^4+5272323840 t^3 x_i^2-47933596800 t^3 x_i\\&-&12312919440 t^3+25284879360 t^2 x_i^2+12642439680 t^2 x_i+1580304960 t^2\Big\}, \\
  T_i &=& \frac{1}{\left(4 t^2-4 t+4 x_i+1\right)^3}\Big\{64 B(t) E(t) t^6-192 B(t) E(t) t^5+192 B(t) E(t) t^4 x_i+240 B(t) E(t) t^4\\&-&384 B(t) E(t) t^3 x_i-160 B(t) E(t) t^3+192 B(t) E(t) t^2 x_i^2+288 B(t)
   E(t) t^2 x_i+60 B(t) E(t) t^2\\&-&192 B(t) E(t) t x_i^2-96 B(t) E(t) t x_i-12 B(t) E(t) t+64 B(t) E(t) x_i^3+48 B(t) E(t) x_i^2\\&+&12 B(t) E(t) x_i+B(t)
   E(t)+128 C(t) E(t) t^4-256 C(t) E(t) t^3+256 C(t) E(t) t^2 x_i\\&+&192 C_i(t) E(t) t^2-256 C(t) E(t) t
   x_i-64 C(t) E(t) t+128 C(t) E(t) x_i^2+64 C(t) E(t) x_i\\&+&8 C(t) E(t)+192 D(t) E(t) t^2-192 D(t)
   E(t) t+192 D(t) E(t) x_i+48 D(t) E(t)+64 D(t) t^2 y_i\\&-&64 D(t) t y_i+64 D(t) x_i y_i+16 D(t) y_i+256 E^3(t)+256
  E^2(t) y_i\Big\}. \\
\end{eqnarray*}}
Using MAGMA \cite{MAGMA}, the specialization $t=\frac{3}{4}$ shows that the rational points $q_i$, $i=0,1,2,3,4,5$, on the elliptic curve $E_t^*$ are independent. According to Silverman's Specialization Theorem, the curve $E_t^*$ has rank $\geq 6$. Therefore, the above procedure yields an infinite family of elliptic curves with a $2$-torsion point and rank at least $6$.

\section{Second construction}
In this section, given a $6$-term sequence of consecutive squares, we establish the existence of infinitely many elliptic curves $C$ described by $y^2=ax^4+bx^2+c$ over $\Q$ with this sequence making up the $x$-coordinates of rational points on $C$.

In fact, given $t\in\Q$ such that $((t-1)^2,d),(t^2,e)$, and $((t+1)^2,f)$ are rational points in $C(\mathbb{Q})$, where $C:y^2=ax^4+bx^2+c$, one sees that
\begin{eqnarray*}
d^2&=&a(t-1)^8+b(t-1)^4+c\\e^2&=&at^8+bt^4+c\\f^2&=&a(t+1)^8+b(t+1)^4+c.
\end{eqnarray*}
The values of $a,b,c\in \Q[d,e,f](t)$ can be determined by solving the above system of linear equations.

Furthermore, if $((t-2)^2,g)$ is a fourth rational point in $C(\Q)$, then the values of $a,b,c$ yield that
{\footnotesize\begin{eqnarray}
\label{equ1}
g^2&=&\frac{3  \left(8 t^8-20 t^7+36 t^6-24 t^5+2 t^4+15 t^3+9 t^2-16 t-10\right)d^2}{t \left(8 t^7+4 t^6+8 t^5+4 t^4+2 t^3+t^2+2 t+1\right)}\nonumber\\& &-\frac{3
   \left(4 t^3-18 t^2+28 t-15\right) \left(2 t^4-2 t^3+7 t^2-2 t+5\right)e^2}{8 t^7+4 t^6+8 t^5+4 t^4+2 t^3+t^2+2 t+1}\nonumber\\& &+\frac{ \left(4 t^3-18 t^2+28
   t-15\right) \left(2 t^5-8 t^4+15 t^3-15 t^2+8 t-2\right)f^2}{t \left(8 t^7+4 t^6+8 t^5+4 t^4+2 t^3+t^2+2 t+1\right)}.\nonumber\\
\end{eqnarray}}
In the quadratic equation above, since $(d,e,f,g)=(1,1,1,1)$ is a solution, one may find a parametric solution $(d,e,f,g)\in \Q[t][p,q,w]$.

Now in order for the elliptic curve $C:y^2=ax^4+bx^2+c$ to possess a $5$-term sequence of consecutive squares we assume the existence of the fifth rational point $((t+2)^2,h)$ in $C(\Q)$. This will yield that
\begin{eqnarray}
\label{eq3}
h^2=Ap^4+Bp^3+Cp^2+Dp+E,\textrm{\hskip20pt} A,B,C,D,E\in \Q[t][q,w].
\end{eqnarray}

The specialization $$q=\frac{8 \left(20 t^6 w-54 t^5 w+130 t^4 w-219 t^3 w+257 t^2 w-132 t w-2 w\right)}{3 \left(48 t^6-96 t^5+420 t^4-496 t^3+756 t^2-312 t+5\right)}$$ kills the discriminant of the the binary quartic (\ref{eq3}). In other words, the algebraic curve defined by eq (\ref{eq3}) is singular. More precisely, a parametric solution $(p,h)$ can be given by
{\footnotesize\begin{eqnarray*}\label{equ3}
h&=&\frac{(t-1) (2 t+1) \left(t^2-2 t+2\right) \left(2 t^2+2 t+1\right)}{3 (2 t-1)^2 \left(2 t^2-2 t+5\right) \left(12 t^3-6 t^2+60 t-1\right)^2}\Big(9 (2 t-1)^2 \left(2 t^2-2 t+5\right)^2 \\&\times&\left(12 t^3-6 t^2+60 t-1\right)^2 p^2-96  t (2 t-1) \left(t^2+4\right) \left(2 t^2-2 t+5\right) \left(4 t^3-42
   t^2+4 t-31\right)\\&\times& \left(12 t^3-6 t^2+60 t-1\right)p w-(2 t-3) \left(2 t^2-6 t+5\right) \Big(1088 t^9+22944 t^8+13680 t^7+179048 t^6\\&+&67104 t^5+400204
   t^4+110908 t^3+211754 t^2-2140 t-15\Big) w^2\Big).
\end{eqnarray*}}
We summarize our findings in the following proposition.
\begin{Proposition}
\label{prop1}
Given a nontrivial sequence of consecutive rational squares $(t+i)^2$, $i=0,\pm 1,\pm 2$, there exist infinitely many elliptic curves of the form $C: y^2 = a(t,p,w)x^4 + b(t,p,w)x^2 + c(t,p,w)$
such that $(t + i)^2$ is the $x$-coordinate of a rational point on $C$.
\end{Proposition}

Now we are looking for elliptic curves containing $6$-term sequences of consecutive squares. So we let $((t+3)^2,k)$ be a point in $C(\Q)$ where $C$ is given by $y^2=ax^4+bx^2+c$. Therefore, one has
\begin{eqnarray}\label{eq:ellipticsurface12}
                   k^2=\overline{A} p^4+\overline{B} p^3w +\overline{C} p^2w^2 +\overline{D} pw^3 + \overline{E}w^4 \end{eqnarray}
                   where the description of $\overline{A}, \overline{B},\overline{C},\overline{D},\overline{E}\in\Q(t)$ can be found, for instance, using \Magma.
\begin{Theorem}
\label{thm2}
The curve $\mathcal{C}:k^2=\overline{A} p^4+\overline{B} p^3 +\overline{C} p^2 +\overline{D} p+ \overline{E}$ defined over $\mathbb{Q}(t)$ is birationally equivalent over $\mathbb{Q}(t)$ to an elliptic curve $\mathcal{E}$  with Mordell-Weil rank $\rank\,\mathcal{E}(\mathbb{Q}(t))\ge 1$.
\end{Theorem}
\begin{Proof}
 The curve $\mathcal C$ is nonsingular as the discriminant is nonzero. Moreover, since $\overline{A}$ is a square in $\Q(t)$, one knows that $\mathcal{C}$ is an elliptic curve over $\Q(t)$. In fact, $\mathcal C$ is birationally equivalent to $\mathcal E$ defined by $y^2=x^3-27Ix-27J$ where $I=12\overline{A}\,\overline{E}-3\overline{B}\,\overline{D}+\overline{C}^2$ and $J=72 \overline{A}\,\overline{C}\,\overline{E} + 9 \overline{B}\,\overline{C}\,\overline{D} -27 \overline{A}\,\overline{D}^2 -27\overline{B}^2\overline{E} -2 \overline{C}^3$ with $P=\displaystyle \left( 3\frac{3\overline{B}^2 -8\overline{A}\,\overline{C}}{4\overline{A}}, 27\frac{\overline{B}^3 + 8\overline{A}^2\overline{D} - 4\overline{A}\,\overline{B}\,\overline{C}}{8\overline{A}^{3/2}}\right)$ in $\mathcal{E}(\mathbb{Q}(t))$, see \cite{StollCremona}.

 One considers the specialization $t=3$ in order to obtain the specialized point {\footnotesize $$\widetilde{P}=\left(\frac{19558022787408000000}{201601},\frac{86476754780118743040000000000}{90518849}\right)$$} of the point $P$ on the specialized elliptic curve
{\footnotesize\begin{eqnarray*}\widetilde{\mathcal{E}}: y^2&=&x^3-\frac{156217789162987774532352000000000000}{40642963201}x\\& &+\frac{22789637573454810302335707893243904000000000000000000}{8193662024284801}.\end{eqnarray*}}
 Using $\mathrm{Magma}$ \cite{MAGMA}, the point $\widetilde{P}$ is a point of infinite order on $\widetilde{\mathcal{E}}$. Therefore, according to Silverman's specialization Theorem, the point $P$ is of infinite order on $\mathcal{E}$.
\end{Proof}
\begin{Theorem} \label{thm1}
 For any $6$-term nontrivial sequence of consecutive squares $(t+i)^2$, $i=0,\pm 1,\pm 2,3$, $t\in\Q$, there is an infinite family of elliptic curves described by $C:y^2=a x^4+bx^2+c$, $a,b,c\in\Q$, such that $(t+i)^2$ is the $x$-coordinate of a rational point in $C(\Q)$. In particular, there are infinitely many elliptic curves described by $y^2=a x^4+bx^2+c$ with $6$-term sequences of consecutive squares.
\end{Theorem}
\begin{Proof}
Fix $t=t_0\in\Q$. The values of $a$, $b$, $c$ in $\Q[p,q,w]$ guarantees the existence of the points $((t-1)^2, d),(t^2,e), ((t+1)^2,f), ((t-2)^2,g)$ in $C(\Q)$, where $d,e,f,g\in\Q[p,q,w]$. Choosing {\footnotesize$$\displaystyle q=\frac{8 \left(20 t^6 w-54 t^5 w+130 t^4 w-219 t^3 w+257 t^2 w-132 t w-2 w\right)}{3 \left(48 t^6-96 t^5+420 t^4-496 t^3+756 t^2-312 t+5\right)}$$} yields the existence of a rational point $((t_0+2)^2,h)$ in $C(\Q)$, see Proposition \ref{prop1}. Theorem \ref{thm2} establishes the existence of infinitely many projective pairs $(p:w)$ for which $((t_0+3)^2,k)$ lies in $C(\Q)$.
\end{Proof}
\hskip-12pt\emph{\bf{Acknowledgements.}}
We would like to thank the referees for many comments, corrections, and
suggestions that helped the authors improve the manuscript.

\end{document}